\documentclass{amsart}

\usepackage{amsfonts, amssymb, amsmath, eucal, verbatim, amsthm, amscd, enumerate}
\usepackage{graphicx}

\newtheorem{theorem}{Theorem}
\newtheorem{lemma}[theorem]{Lemma}
\newtheorem{proposition}[theorem]{Proposition}

\newtheorem{corollary}[theorem]{Corollary}

\theoremstyle{definition}

\theoremstyle{remark}
\newtheorem*{remark}{Remark}

\newtheorem*{acknowledgement}{Acknowledgement}

\def\R{{\mathbb R}}

\newcommand{\op}[1]{\operatorname{#1}}

\title{Some sharp bilinear space-time estimates for the wave equation}
\author{Neal Bez} 
\address{Neal Bez, Department of Mathematics, Graduate School of Science and Engineering,
Saitama University, Saitama 338-8570, Japan}
\email{nealbez@mail.saitama-u.ac.jp}

\author{Chris Jeavons}
\address{Chris Jeavons, School of Mathematics, University of Birmingham, Edgbaston, Birmingham, B15 2TT, UK}
\email{jeavonsc@maths.bham.ac.uk}

\author{Tohru Ozawa}
\address{Tohru Ozawa, Department of Applied Physics, Waseda University, Tokyo 169-8555, Japan}
\email{txozawa@waseda.jp}


\date{\today}

\begin{document}

\begin{abstract}
We prove a family of sharp bilinear space-time estimates for the half-wave propagator $e^{it\sqrt{-\Delta}}$. As a consequence, for radially symmetric initial data, we establish sharp estimates of this kind for a range of exponents beyond the classical range. 
\end{abstract}
\maketitle
\section{Introduction and statement of results}
For $d\geq 2$, suppose that $u,v$ satisfy the homogeneous wave equations 
\[
\square u = \square v = 0,
\]
where $\square = \Delta_x - \partial_t^2$ with $(t,x) \in \mathbb{R} \times \mathbb{R}^{d}$, with initial conditions
\[
(u(0),\partial_tu(0)) = (u_0,u_1), \qquad  (v(0),\partial_tv(0)) = (v_0,v_1).
\]
We consider the bilinear global space-time estimates
\begin{equation}\label{FK}
\|D^{\beta_0}D_-^{\beta_-}D_+^{\beta_+}(uv)\|_{L^2(\R^{d+1})}\leq C \|(u_0,u_1)\|_{\dot{H}^{\alpha_1}\times\dot{H}^{\alpha_1-1}}\|(v_0,v_1)\|_{\dot{H}^{\alpha_2}\times\dot{H}^{\alpha_2-1}},
\end{equation}
where $\dot{H}^s(\R^d)$ denotes the usual homogeneous Sobolev space with norm
\[\|f\|_{\dot{H}^s(\R^d)}:=\|D^sf\|_{L^2(\R^d)}.\]
Here, we adopt the standard notation $D:=\sqrt{-\Delta_x}$, and the operators $D_+$ and $D_-$ are defined as
\[
\widetilde{D_{\pm}f}(\tau,\xi):=\big||\tau|\pm|\xi|\big|\widetilde{f}(\tau,\xi)
\] 
using the space-time Fourier transform
\[
\widetilde{f}(\tau,\xi):=\int_{\R^{d+1}}e^{-it\tau-ix\cdot\xi}f(t,x)\,\mathrm{d}x\mathrm{d}t.
\] 
Estimates \eqref{FK} go back to work of Beals in \cite{Beals} and Klainerman--Machedon in \cite{KM'93}; certain special cases were established and used in the study of the null forms for some nonlinear wave equations. Further special cases were proved, and some applications presented, in \cite{KM'95}, \cite{KM'96}, \cite{KM'96-2}, \cite{KM'97}, \cite{KM'97-2} and \cite{KS} among others, and the full range of exponents $(\beta_0,\beta_-,\beta_+,\alpha_1,\alpha_2)$ for which \eqref{FK} holds was found by Foschi--Klainerman in \cite{FK}. 

Although not the focus of the present paper, we mention that generalisations of \eqref{FK} where $L^2$ norm on the left-hand side is replaced by a mixed space-time $L_t^qL_x^r$ norm have also been extensively studied. The first progress for $q = r < 2$ was obtained by Bourgain \cite{Bourgain2}, considering bilinear estimates without the multiplier weights but with separated frequency supports, and the sharp estimates (in the sense of exponents) were later obtained by Wolff \cite{Wolff} and subsequently Tao \cite{Tao} (generalised by Lee \cite{Lee}). With the multipler weights, the full range of exponents in the mixed-norm case was obtained (up to endpoints) by Lee--Vargas \cite{LV} for $d \geq 4$, and by Lee--Rogers--Vargas \cite{LRV} when $d=3$, leaving only some non-endpoint cases when $d=2$.


If $\square u = 0$ then we have the standard decomposition into ($+$) and ($-$) waves $u=e^{itD}f_++e^{-itD}f_-$, where $f_\pm$ are given by
\begin{equation}\label{upm}
u(0)=f_++f_-, \quad \partial_t u(0)=iD(f_+-f_-),
\end{equation}
and the half-wave propagator $e^{itD}$ is defined on the Schwartz space using (spatial) Fourier inversion to be
\[
e^{itD}f(x):=\frac{1}{(2\pi)^d}\int_{\R^d}e^{ix\cdot\xi}e^{it|\xi|}\widehat{f}(\xi)\,\mathrm{d}\xi, \qquad (t,x)\in\R\times\R^d.
\]
Using this decomposition, the inequalities \eqref{FK} can be reduced to the following bilinear estimates for the propagator $e^{itD}$:
\begin{equation}\label{FK1s}
\|D^{\beta_0}D_-^{\beta_-}D_+^{\beta_+}\left(e^{itD}f e^{itD}g\right)\|_{L^2(\R^{d+1})}\leq C \|f\|_{\dot{H}^{\alpha_1}(\R^d)}\|g\|_{\dot{H}^{\alpha_2}(\R^d)},
\end{equation}
and
\begin{equation}\label{FK1sbar}
\|D^{\beta_0}D_-^{\beta_-}D_+^{\beta_+}\left(e^{itD}f \overline{e^{itD}g}\right)\|_{L^2(\R^{d+1})}\leq C \|f\|_{\dot{H}^{\alpha_1}(\R^d)}\|g\|_{\dot{H}^{\alpha_2}(\R^d)}.
\end{equation}
In addition to the works cited above, estimates related to \eqref{FK1s} appear in the book \cite{Bourgain}, and applications to well-posedness problems for some nonlinear wave equations are also discussed there.

Closely related to the estimates \eqref{FK1s} and \eqref{FK1sbar}, at least in the symmetric cases $\beta_-=\beta_+$ and $\alpha_1=\alpha_2$, is the following inequality:
\begin{align}\label{BR}
& \|e^{itD} f e^{itD} g \|_{L^2(\R^{d+1})}^2 \nonumber \\
& \leq    \frac{2^{\frac{5-7d}{2}} \pi^{\frac{2-5d}{2}}} {\Gamma(\frac{d}{2})}   \int_{\mathbb{R}^{2d}}|\widehat{f}(y_1)|^2|\widehat{g}(y_2)|^2 (|y_1||y_2|)^{\frac{d-1}{2}} (1-y_1^\prime\cdot y_2^\prime)^{\frac{d-3}{2}}\,\mathrm{d}y_1\mathrm{d}y_2,
\end{align}
where $y' := |y|^{-1}y \in \mathbb{S}^{d-1}$. When $d\geq 3$ the constant is sharp, and a full characterisation of extremising initial data is known; this was first proved by Foschi \cite{Foschi} in the case\footnote{This is stated for the case $f=g$, but the same argument in \cite{Foschi} allows $f$ and $g$ to be different.} $d=3$, and for arbitrary dimensions in \cite{BR}. The question of determining sharp constants and identifying extremisers for weighted linear and bilinear space-time estimates for dispersive and wave-like propagators has been studied in a number of recent papers. In addition to \cite{BR} and \cite{Foschi} already mentioned, see \cite{BBJP}, \cite{Carneiro}, \cite{HZ} and \cite{OT} (as well as \cite{Foschi}, again) for the Schr\"odinger equation, \cite{Jeavons}, \cite{OR2} and \cite{Quilodran} for the Klein--Gordon equation, \cite{BezJ} for the wave equation, and \cite{BezSugimoto}, \cite{BezSugimoto2} and \cite{OR} as well as further references contained in these papers, for more general propagators and spatial weights in the linear setting on $L^2$.

Our main result is the following one-parameter family of sharp bilinear inequalities in the case of ($+-$) waves. Before proceeding, we introduce some notation. For $\beta\in\R$, we write
\[
\mathcal{I}_\beta(f,g)=\int_{\R^{2d}}|\widehat{f}(y_1)|^2|\widehat{g}(y_2)|^2 (|y_1||y_2|)^{\frac{d-1}{2} + 2\beta} \left(1 -y_1'\cdot y_2'\right)^{\frac{d-3}{2}+2\beta}\,\mathrm{d}y_1\mathrm{d}y_2,
\]
defined for suitable functions $f$ and $g$. Also, we write $|\square| = D_- D_+$, so that
\[
\widetilde{|\square|^\beta u}(\tau,\xi):=|\tau^2-|\xi|^2|^{\beta}\widetilde{u}(\tau,\xi).
\]
Finally, we set $\beta_d := \max\{\frac{1-d}{4}, \frac{2-d}{2} \}$, so that $\beta_d = \frac{1-d}{4}$ when $d \geq 3$, and $\beta_2 = 0$.
\begin{theorem}\label{thm1}
Let $d\geq 2$ and $\beta>\frac{1-d}{4}$. Then
\begin{equation}\label{b-eq}
\||\square|^\beta(e^{itD}f\overline{e^{itD}g})\|_{L^2(\mathbb{R}^{d+1})}^2\leq \mathbf{W}(\beta,d)\mathcal{I}_\beta(f,g)
\end{equation}
holds with constant
\[
\mathbf{W}(\beta,d)= 2^{\frac{1-5d + 4\beta}{2}} \pi^{\frac{1-5d}{2}}  \frac{\Gamma( \frac{d-1}{2}  + 2\beta)}{\Gamma(d-1+2\beta)}.
\]
For any $\beta > \beta_d$, the constant is optimal and equality holds if and only if
\begin{equation}\label{maxdef}
\widehat{f}(\xi)=\lambda\widehat{g}(\xi)=\frac{e^{a|\xi|+b\cdot\xi+c}}{|\xi|}
\end{equation}
for some $\lambda, a, c\in\mathbb{C}$ with $\operatorname{Re}a<0$, and $b\in\mathbb{C}^d$ with $|\operatorname{Re}b|<-\operatorname{Re}a$.
\end{theorem}
Estimate \eqref{b-eq} implicitly holds for initial data such that $\mathcal{I}_\beta(f,g)$ is finite, and the necessity of the lower bound $\beta > \frac{1-d}{4}$ is evident from the expression for $\mathbf{W}(\beta,d)$. In order to state that the constant is sharp, the additional lower bound $\beta > \frac{2-d}{2}$ arises by evaluating $\mathcal{I}_\beta(f,g)$ on initial data satisfying \eqref{maxdef}.

When $\beta=0$, inequality \eqref{b-eq} recovers \eqref{BR} since the duplication formula for the gamma function gives
\[
\mathbf{W}(0,d) = 2^{\frac{1-5d}{2}} \pi^{\frac{1-5d}{2}}  \frac{\Gamma(\frac{d-1}{2})}{\Gamma(d-1)}  =  \frac{2^{\frac{5-7d}{2}} \pi^{\frac{2-5d}{2}}} {\Gamma(\frac{d}{2})},
\]
and the $(++)$ waves on the left-hand side can be freely switched with $(+-)$ waves. The case $\beta=\frac{3-d}{4}$ in Theorem \ref{thm1} is obviously special since the power of the angular weight defining $\mathcal{I}_\frac{3-d}{4}$ is zero and so by Plancherel's theorem,
\begin{equation} \label{e:threshold}
\mathcal{I}_{\frac{3-d}{4}}(f,g)=(2\pi)^{2d}\|f\|_{\dot{H}^{\frac{1}{2}}(\mathbb{R}^d)}^2\|g\|_{\dot{H}^{\frac{1}{2}}(\mathbb{R}^d)}^2.
\end{equation}
Therefore, \eqref{b-eq} recovers inequality \eqref{FK1sbar} in the case 
$$
(\beta_0,\beta_-,\beta_+,\alpha_1,\alpha_2)=\left(0,\frac{3-d}{4},\frac{3-d}{4},\frac{1}{2},\frac{1}{2}\right)
$$
with sharp constant. In this sense Theorem \ref{thm1} unifies \eqref{BR} with certain cases of \eqref{FK1sbar}. Of course, one may also ask for a similar unification of \eqref{BR} with certain cases of \eqref{FK1s}, by establishing analogous sharp estimates to \eqref{b-eq} for $(++)$ waves. It turns out that such estimates are already essentially contained in \cite{BR} and the argument given there for $\beta = 0$ easily extends to general $\beta$. It is not the case that \eqref{b-eq} follows from the arguments in \cite{BR} and hence the main focus in this paper is the more difficult $(+-)$ case. For completeness, we provide a new proof of the analogous estimates to \eqref{b-eq} for $(++)$ waves in the final section, with an argument which clearly elucidates the additional technical difficulty in the $(+-)$ case.

Inequality \eqref{BR} has been shown to imply sharp forms of certain linear $L^4$ Sobolev--Strichartz estimates for the wave equation, by controlling the right-hand side of \eqref{BR} in a sharp manner in terms of classical homogeneous Sobolev norms; this is done in \cite{BR} in five space dimensions, and more recently in \cite{BezJ} in four space dimensions. Here we observe that if the initial data $f$ and $g$ are both radially symmetric, we can use polar coordinates to calculate 
\[
\mathcal{I}_\beta(f,g) = 2^{\frac{7(d-1)}{2} + 2\beta} \pi^{\frac{4d-1}{2}} \frac{\Gamma(\frac{d}{2}) \Gamma(d-2+2\beta)}{\Gamma(\frac{3d-5}{2} + 2\beta)} \|f\|_{\dot{H}^{\frac{d-1}{4} + \beta}(\R^d)}^2\|g\|_{\dot{H}^{\frac{d-1}{4} + \beta}(\R^d)}^2
\] 
and the following sharp estimates valid in \emph{all} dimensions may be obtained immediately from Theorem \ref{thm1}.
\begin{corollary} \label{c:radial}
Let $d \geq 2$, $\beta > \beta_d$ and let
\begin{equation} \label{e:Cdefn}
\mathbf{C}(\beta,d)=\frac{2^{d-3+4\beta}\Gamma\left(\frac{d}{2}\right)\Gamma\left(\frac{d-1}{2} + 2\beta\right)}{\pi^{\frac{d}{2}}(d-2+2\beta)\Gamma\left(\frac{3d-5}{2}+2\beta\right)}.
\end{equation}
If $f$ and $g$ are radially symmetric, then 
\begin{equation} \label{e:radialall}
\left\||\square|^\beta (e^{itD}f\overline{e^{itD}g})\right\|_{L^2(\R^{d+1})}^2\leq \mathbf{C}(\beta,d)\|f\|_{\dot{H}^{\frac{d-1}{4} + \beta}(\R^d)}^2\|g\|_{\dot{H}^{\frac{d-1}{4} + \beta}(\R^d)}^2,
\end{equation}
where the constant is sharp and equality holds if and only if $f=\lambda g$ satisfy \eqref{maxdef} with $\lambda, c\in\mathbb{C}$, and $a \in \mathbb{C}$ such that $\operatorname{Re}a < 0$, and $b=0$.
\end{corollary}
Without the additional hypothesis of radial symmetry, the estimate \eqref{e:radialall} fails to hold for \emph{any} finite constant when $\beta < \frac{3-d}{4}$.
\begin{proposition}\label{p:counter}
If $\beta<\frac{3-d}{4}$, then for any $A>0$ there exists $f\in\dot{H}^{\frac{d-1}{4}+\beta}(\R^d)$ such that
\[\frac{\left\||\square|^\beta|e^{itD}f|^2\right\|_{L^2(\R^{d+1})}}{\|f\|_{\dot{H}^{\frac{d-1}{4}+\beta}(\R^d)}^2}>A.\]
\end{proposition}
Thus, \eqref{e:radialall} shows that the range of admissible exponents is widened when the initial data are radially symmetric. Moreover, our expression for the sharp constant $\mathbf{C}(\beta,d)$ clearly indicates that the range $\beta > \beta_d$ cannot be further widened, even under the radial hypothesis on the initial data. Results of this type for initial data restricted to the radial case for estimates similar to \eqref{FK} were proved by Klainerman--Machedon in \cite{KM'93}, and for the estimates \eqref{FK} themselves by Foschi in \cite{Foschi3} (in the case $\beta_- = \beta_+  = 0$). The fact that certain estimates, including Strichartz-type estimates, improve on radially symmetric data is well-known and has been studied more widely; see for example \cite{BGO}, \cite{CO}, \cite{FW}, \cite{Hidano}, \cite{HK}, \cite{Sogge}, \cite{Sterbenz}. 

Regarding Proposition \ref{p:counter}, we remark that in Example 5.1 of \cite{FK} it is shown that $\beta_- \geq \frac{3-d}{4}$ is a necessary condition for \eqref{FK} to hold by showing that \eqref{FK1s} fails in the ($++$) case. To prove Propositon \ref{p:counter} we must consider the ($+-$) case, and we use a different argument (using initial data of the form \eqref{maxdef}) to establish this.

For a restricted range of $\beta$ we are able to extend Corollary \ref{c:radial} to general initial data without the assumption of radial symmetry, by applying some of the techniques from \cite{BezJ} to make a bound on the right-hand side of \eqref{b-eq} which does not incur an overall loss of sharpness in the constant. In order to state this, we introduce a collection of maps $T_\beta:\dot{H}^{\frac{d-1}{4}+\beta}(\R^d)\rightarrow L^1(\mathbb{S}^{d-1})$ defined by
\[T_\beta f(\omega)=\frac{1}{(2\pi)^d}\int_0^\infty |\widehat{f}(r\omega)|^2 r^{\frac{3d-3}{2}+2\beta}\,\mathrm{d}r,\]
for $\omega\in\mathbb{S}^{d-1}$. 
\begin{corollary}\label{cor4}
Let $d \geq 2$, $\beta \in (\beta_d, \tfrac{5-d}{4}]$ and let $\mathbf{C}(\beta,d)$ be given by \eqref{e:Cdefn}.

\emph{(i).} If $\beta>\frac{3-d}{4}$, then 
\begin{equation}\label{e:alpha}
\left\||\square|^\beta|e^{itD}f|^2\right\|_{L^2(\R^{d+1})}^2\leq \mathbf{C}(\beta,d)\|f\|_{\dot{H}^{\frac{d-1}{4} + \beta}(\R^d)}^4
\end{equation}
and the constant is sharp, with equality if and only if $f$ satisfies \eqref{maxdef} with $c\in\mathbb{C}$, and $a \in \mathbb{C}, b \in \mathbb{C}^d$ are such that $\operatorname{Re}a < 0$ and $\operatorname{Re}b=0$.

\emph{(ii).} If $\beta\leq\frac{3-d}{4}$ and $p=\frac{2(d-1)}{3d-5+4\beta}$, then
\begin{equation}\label{spec}
\left\||\square|^\beta (e^{itD}f\overline{e^{itD}g})\right\|_{L^2(\R^{d+1})}^2\leq \mathbf{C}(\beta,d)|\mathbb{S}^{d-1}|^{\frac{3-d-4\beta}{d-1}}\|T_\beta f\|_{L^p}\|T_\beta g\|_{L^p}
\end{equation}
and the constant is sharp, with equality if and only if $f=\lambda g$ satisfy \eqref{maxdef} with $\lambda, a, c\in\mathbb{C}$ and $b\in\mathbb{C}^d$ where $|\operatorname{Re}b|<-\operatorname{Re}a$. 
\end{corollary}
Corollary \ref{cor4} extends Corollary \ref{c:radial} in two ways. Estimate \eqref{e:alpha} extends \eqref{e:radialall} to general initial data in $\dot{H}^{\frac{d-1}{4}+\beta}(\R^d)$ for $\beta \in (\frac{3-d}{4},\frac{5-d}{4}]$ (and in the symmetric case $f=g$). The upper threshold $\frac{5-d}{4}$ is a consequence of our method and it is quite possible the sharp estimate can be extended beyond this threshold (however, it seems that a new key idea is necessary for this). 

For $\beta \in (\beta_d,\frac{3-d}{4})$, the estimate \eqref{spec} provides a natural substitute for the failure of \eqref{e:alpha} which is valid for all initial data in $\dot{H}^{\frac{d-1}{4}+\beta}(\R^d)$, thus extending \eqref{e:radialall} for such $\beta$. Indeed, for such $\beta$ we have $p=\frac{2(d-1)}{3d-5+4\beta} > 1$ and, by Plancherel's theorem and H\"older's inequality,
\begin{equation}\label{PH}
\|f\|_{\dot{H}^{\frac{d-1}{4}+\beta}(\R^d)}^2\leq |\mathbb{S}^{d-1}|^{\frac{3-d-4\beta}{2(d-1)}}\|T_\beta f\|_{L^p(\mathbb{S}^{d-1})},
\end{equation} 
with equality when $f$ is radially symmetric.

\emph{Organisation.} In Section \ref{s:cor} we prove Proposition \ref{p:counter} and Corollary \ref{cor4}, and in Section \ref{s:proofs} we prove Theorem \ref{thm1} for $(+-)$ waves. Finally, in Section 4, we present analogous results for $(++)$ waves.

\section{Proofs of Proposition \ref{p:counter} and Corollary \ref{cor4}}\label{s:cor}
We begin this section with the proof of Corollary \ref{cor4} following arguments in \cite{BezJ}. Before proceeding, we introduce notation
\[
H_\lambda(f,g):=\int_{\mathbb{S}^{d-1}} \int_{\mathbb{S}^{d-1}}f(\omega_1)\overline{g(\omega_2)}|\omega_1-\omega_2|^{-\lambda}\,\mathrm{d}\omega_1\mathrm{d}\omega_2
\]
where $|\cdot|$ in this context means chordal distance on $\mathbb{S}^{d-1}$ (that is, euclidean distance on $\R^d$) and we fix $\lambda=3-d-4\beta$. Of course, $\lambda \geq 0$ if and only if $\beta \leq \frac{3-d}{4}$, and it is reasonable to expect the behaviour of $H_\lambda$ to change according to the sign of $\lambda$.
\proof[Proof of Corollary \ref{cor4}.]
The proof rests on the observation
\begin{equation} \label{HLStrick}
\mathcal{I}_\beta(f,g)=\frac{(2\pi)^{2d}}{2^{\frac{d-3}{2}+2\beta}}H_\lambda(T_\beta f, T_\beta g)
\end{equation}
which is easily verified.

We begin with the case $\beta \in (\frac{3-d}{4},\frac{5-d}{4}]$ for each $d \geq 2$, which means $\lambda \in [-2,0)$. In this case, we use the following result from \cite{BezJ} to bound $H_\lambda$.
\begin{lemma}\label{claim1}
Let $d\geq 2$, $\lambda \in [-2,0)$, and let $g$ be any $L^1$ function on $\mathbb{S}^{d-1}$. Then,
\[
H_\lambda(g,g) \leq 2^{2d-5+4\beta} \pi^{-\frac{1}{2}} \frac{\Gamma(d-2+2\beta) \Gamma(\frac{d}{2})}{\Gamma(\frac{3d-5}{2} + 2\beta)}   \left|\int_{\mathbb{S}^{d-1}}g\right|^2,
\]
with equality if $g$ is constant. If, in addition $\lambda>-2$, then equality holds only if $g$ is constant.
\end{lemma} 
\begin{remark}
Strictly speaking, Lemma \ref{claim1} was proved in \cite{BezJ} only for $d\geq 3$, however the same argument also proves it for $d=2$, using the standard convention $|\mathbb{S}^0|=2$.  
\end{remark}
Since $T_\beta f\geq 0$ and we have that
\[
\|T_\beta f\|_{L^1(\mathbb{S}^{d-1})}=\|f\|_{\dot{H}^{\frac{d-1}{4}+\beta}(\R^d)}^2
\]
by Plancherel's theorem. Thus, Lemma \ref{claim1} implies that the right-hand side of \eqref{b-eq} with $f=g$ is at most
\begin{equation*}
\mathbf{C}(\beta,d)  \|f\|_{\dot{H}^{\frac{d-1}{4} + \beta}(\R^d)}^4,
\end{equation*}
with equality if $T_\beta f$ is constant on $\mathbb{S}^{d-1}$, which certainly happens when $|\widehat{f}|$ is radial. An example of such a function of the form \eqref{maxdef} is given by $\widehat{f}=\frac{e^{-|\cdot|}}{|\cdot|}$ and this establishes sharpness of the constant in \eqref{e:alpha}. 

The characterisation of extremisers for inequality \eqref{e:alpha} follows as in \cite{BezJ} without new arguments, so we give only brief details here. The idea is that any extremiser $f$ for \eqref{e:alpha} must be of the form \eqref{maxdef} since we applied \eqref{b-eq} to obtain \eqref{e:alpha}, and for such $f$ we have
\begin{align*}
T_\beta f(\omega) = e^{2\op{Re}c}\int_0^\infty e^{r(2\op{Re}a+2\op{Re}b\cdot\omega)}r^{\frac{3d-7}{2}+2\beta}\,\mathrm{d}r. 
\end{align*}
Hence
\begin{equation} \label{e:Tbeta}
T_\beta f(\omega) =\frac{C}{\left(1+\xi\cdot\omega\right)^{\frac{3d-5}{2}+2\beta}} 
\end{equation}
for some positive constant $C$ and where $\xi :=\frac{\op{Re}b}{\op{Re}a}\in\R^d$ satisfies $|\xi|<1$ by the assumptions on $a$ and $b$. For $\beta \in (\frac{3-d}{4},\frac{5-d}{4})$ (so that $\lambda \in (-2,0)$), the characterisation of extremisers in Lemma \ref{claim1} means $T_\beta f$ is constant and so $\xi$ must be zero. Hence $\op{Re}b = 0$ and this gives the claimed characterisation of extremisers for \eqref{e:alpha}. The remaining case $\beta = \frac{5-d}{4}$ (so that $\lambda = -2$) is more subtle because there are non-constant extremisers for the estimate in Lemma \ref{claim1}. However, the argument in Section 3 of \cite{BezJ} yields the same conclusion\footnote{Strictly speaking the argument in Section 3 of \cite{BezJ} is for the case $d=5$ but it can easily be shown to generalise to all $d \geq 2$.} that $\op{Re}b = 0$, and we omit the details.

Now we move to the case $\beta \in (\beta_d,\frac{3-d}{4}]$. When $\beta = \frac{3-d}{4}$ (or $\lambda = 0$), \eqref{spec} and the characterisation of extremisers follows from Theorem \ref{thm1} and \eqref{e:threshold}.

In the case $\beta \in (\beta_d,\frac{3-d}{4})$, we have $\lambda \in (0,d-1)$ so we apply the sharp Hardy--Littlewood--Sobolev inequality on the sphere and characterisation of extremisers (due to Lieb \cite{Lieb}) to obtain
\begin{equation}\label{HLS-ineq}
H_\lambda(T_\beta f,T_\beta g)\leq \pi^{\frac{\lambda}{2}}\frac{\Gamma\left(\tfrac{d-1-\lambda}{2}\right)}{\Gamma\left(d-1-\tfrac{\lambda}{2}\right)}\left(\frac{\Gamma\left(d-1\right)}{\Gamma\left(\tfrac{d-1}{2}\right)}\right)^{1-\frac{\lambda}{d-1}}\|T_\beta f\|_{L^p}\|T_\beta g\|_{L^p},
\end{equation}
for 
\[p:=\frac{2(d-1)}{2(d-1)-\lambda}=\frac{2(d-1)}{3d-5+4\beta}.\]
Equality holds in \eqref{HLS-ineq} if and only if there exist $C_0, C_1\in\mathbb{C}$ and $\xi\in\R^d$ with $|\xi|<1$ such that
\begin{equation}\label{HLS-ext}
T_\beta f(\omega)=\frac{C_0}{\left(1+\xi\cdot\omega\right)^{\frac{2(d-1)-\lambda}{2}}},\qquad T_\beta g(\omega)=\frac{C_1}{\left(1+\xi\cdot\omega\right)^{\frac{2(d-1)-\lambda}{2}}}.
\end{equation}
That inequality \eqref{spec} is true with the stated constant now follows immediately from estimate \eqref{b-eq} in Theorem \ref{thm1} and using \eqref{HLS-ineq}. The optimality of the constant in \eqref{spec} and characterisation of extremisers is deduced from the additional observation that the functions $f, g$ which exhaust the extremals for \eqref{b-eq}, given by \eqref{maxdef}, also satisfy \eqref{HLS-ext}; this follows immediately from \eqref{e:Tbeta}.\endproof

\proof[Proof of Proposition \ref{p:counter}] Fix $f=f_{\delta}$ such that $\widehat{f}$ is a function of the form \eqref{maxdef} with $c=0$, $a=-1$ and $b=(1-\delta)\mathbf{e}_1$ with $0 < \delta < \frac{1}{100}$; that is 
\begin{equation*}
\widehat{f}(\xi)=\frac{e^{-|\xi|+(1-\delta)\mathbf{e}_1 \cdot\xi}}{|\xi|}.
\end{equation*}
As usual $\mathbf{e}_1 = (1,0,\ldots,0)$ denotes the first basis vector in $\R^d$. For the rest of this section we shall let $C$ denote an arbitrary positive constant which may depend on $d$ and $\beta$ but not on $\delta$, and we use $x\lesssim y$ (respectively, $x\gtrsim y$) to mean $x\leq C y$ ($x\geq C y$), where $C$ may be different even in a single chain of inequalities. We denote by $x \sim y$ if both $x \lesssim y$ and $x \gtrsim y$ hold.

From \eqref{e:Tbeta} we have that 
\[
T_\beta f(\omega)=\frac{C}{\left(1-(1-\delta)\mathbf{e}_1\cdot\omega\right)^{\frac{3d-5}{2}+2\beta}} = \frac{C}{\left(1-(1-\delta)\mathbf{e}_1\cdot\omega\right)^{\frac{d-1}{p}}},
\] 
where $p$ is as in \eqref{spec}; note that $\beta \in (\beta_d,\frac{3-d}{4})$ implies that $p \in (1,2)$. Since $f$ is extremal for \eqref{spec} it is enough to show that
\[
\frac{\|T_\beta f\|_{L^p(\mathbb{S}^{d-1})}}{\|f\|_{\dot{H}^{\frac{d-1}{4}+\beta}(\R^d)}^2}=\frac{\|T_\beta f\|_{L^p(\mathbb{S}^{d-1})}}{\|T_\beta f\|_{{L^1(\mathbb{S}^{d-1})}}}\rightarrow\infty
\]
as $\delta \rightarrow 0$. By changing variables, we have that
\[ 
\|T_\beta f\|_{L^p(\mathbb{S}^{d-1})}=C\left(\int_{-1}^1\frac{(1-t^2)^{\frac{d-3}{2}}}{(1-(1-\delta)t)^{d-1}}\,\mathrm{d}t\right)^{\frac{1}{p}}=:I_1(\delta),
\]
and
\[
\|T_\beta f\|_{L^1(\mathbb{S}^{d-1})}=C\int_{-1}^1\frac{(1-t^2)^{\frac{d-3}{2}}}{(1-(1-\delta)t)^{\frac{d-1}{p}}}\,\mathrm{d}t=:I_2(\delta),
\]
so to complete the proof of Proposition \ref{p:counter} it therefore suffices to prove that
\begin{equation}\label{lb}
I_1(\delta) \gtrsim \delta^{\frac{1-d}{2p}},
\end{equation}
and
\begin{equation}\label{ub}
I_2(\delta) \lesssim \delta^{\frac{1-d}{p}+\frac{d-1}{2}}.
\end{equation}

For this, we let $\sigma > \frac{d-1}{2}$ and change variables ($s = \delta^{-1}(1 - (1-\delta)t$)) to obtain
\begin{align*}
\int_{0}^1\frac{(1-t^2)^{\frac{d-3}{2}}}{(1-(1-\delta)t)^{\sigma}}\,\mathrm{d}t \sim \delta^{\frac{d-1}{2} - \sigma} \int_1^{\frac{1}{\delta}} \frac{(s-1)^{\frac{d-3}{2}}}{s^\sigma} \, \mathrm{d}s  \sim \delta^{\frac{d-1}{2} - \sigma}.
\end{align*}
Applying this with $\sigma = d-1$ immediately yields \eqref{lb}. For \eqref{ub}, by a simple change of variables and the fact that if $t>0$ then $1+(1-\delta)t > 1-(1-\delta)t$, we have
\[
I_2(\delta) \lesssim \int_0^1 \frac{(1-t^2)^{\frac{d-3}{2}}}{(1-(1-\delta)t)^{\frac{d-1}{p}}}\,\mathrm{d}t,
\]
from which \eqref{ub} follows by taking $\sigma = \frac{d-1}{p}$. \endproof

\section{Proof of Theorem \ref{thm1}}\label{s:proofs}
\subsection{Proof of the sharp inequality \eqref{b-eq}}
We broadly follow the strategy in \cite{BBJP} (which in turn was based on the approach in \cite{BBI}) to reduce to evaluating a certain integral over a submanifold of $\R^{2d}$. By Plancherel's theorem and using the relabelling $(x_1,x_2,y_1,y_2)\mapsto(x_1,y_2,y_1,x_2)$,
\begin{align*}
& 2^{-2\beta}(2\pi)^{3d-1}\||\square|^\beta (e^{itD}f\overline{e^{itD}g})\|_{L^2(\R^{d+1})}^2 \\
& = \int_{\R^{4d}}(|y_1||y_2|+y_1\cdot y_2)^{2\beta}\widehat{f}(x_1)\overline{\widehat{g}(-x_2)}\overline{\widehat{f}(y_1)}\widehat{g}(-y_2)\delta\binom{-|x_1|+|x_2|+|y_1|-|y_2|}{x_1+x_2-y_1-y_2}\,\mathrm{d}x\mathrm{d}y \\
& = \int_{\R^{4d}}(|y_1||x_2|+y_1\cdot x_2)^{2\beta}\widehat{f}(x_1)\overline{\widehat{g}(-y_2)}\overline{\widehat{f}(y_1)}\widehat{g}(-x_2)\delta\binom{-|x_1|+|y_2|+|y_1|-|x_2|}{x_1+y_2-y_1-x_2}\,\mathrm{d}x\mathrm{d}y \\
& = \int_{\R^{4d}}(|y_1||x_2|+y_1\cdot x_2)^{2\beta}\frac{\widehat{F}(x)\overline{\widehat{F}(y)}}{\left(|x_1||x_2||y_1||y_2|\right)^{\frac{1}{2}}}\delta\binom{|y_1|+|y_2|-|x_1|-|x_2|}{x_1+y_2-y_1-x_2}\,\mathrm{d}x\mathrm{d}y
\end{align*}
for $\widehat{F}(y):=|y_1|^{\frac{1}{2}}|y_2|^{\frac{1}{2}}\widehat{f}(y_1)\widehat{g}(-y_2)$, for $y=(y_1,y_2)\in\R^{2d}$, and where we have used that if $\tau=|y_1|-|y_2|$ and $\xi=y_1+y_2$ then 
\[
\left|\tau^2-|\xi|^2\right|=2\left|-|y_1||y_2|-y_1\cdot y_2\right|=2(|y_1||y_2|+y_1\cdot y_2).
\]
We now define a non-negative function on $\R^{4d}$,
\[\Psi_{x,y}=\Psi(x,y)=\left(\frac{|x_1||x_2|}{|y_1||y_2|}\right)^{\frac{1}{2}}.\]
Taking real parts and then applying the arithmetic-geometric mean inequality to $\widehat{F}(x)\Psi(x,y)^{\frac{1}{2}}$ and $\widehat{F}(y)\Psi(x,y)^{-\frac{1}{2}}$, it follows that
\begin{align*}
& 2^{-2\beta}(2\pi)^{3d-1}\||\square|^\beta(e^{itD}f\overline{e^{itD}g})\|_{L^2(\R^{d+1})}^2 \\
& \leq \frac{1}{2}\int_{\R^{4d}}(|y_1||x_2|+y_1\cdot x_2)^{2\beta}\frac{|\widehat{F}(y)|^2\Psi_{x,y}^{-1}+|\widehat{F}(x)|^2\Psi_{x,y}}{\left(|x_1||x_2||y_1||y_2|\right)^{\frac{1}{2}}}\delta\binom{|y_1|+|y_2|-|x_1|-|x_2|}{x_1+y_2-y_1-x_2}\,\mathrm{d}x\mathrm{d}y \\
& = \int_{\R^{4d}}(|y_1||x_2|+y_1\cdot x_2)^{2\beta}\frac{|\widehat{F}(y)|^2}{|x_1||x_2|}\delta\binom{|y_1|+|y_2|-|x_1|-|x_2|}{x_1+y_2-y_1-x_2}\,\mathrm{d}x\mathrm{d}y,
\end{align*}
and equality holds if and only if
\[\widehat{F}(x)\Psi(x,y)^{\frac{1}{2}}=\widehat{F}(y)\Psi(x,y)^{-\frac{1}{2}},\]
or equivalently
\begin{equation}\label{f-eq}
|x_1||x_2|\widehat{f}(x_1)\widehat{g}(-x_2)=|y_1||y_2|\widehat{f}(y_1)\widehat{g}(-y_2)
\end{equation}
almost everywhere on the support of the delta measures; for example this equation is satisfied by $f,g$ given by \eqref{maxdef}. 

We now define 
\begin{equation} \label{e:Idefn}
\mathrm{I}_\beta(y):=\int_{\R^d}\int_{\R^d}\frac{(|y_1||x_2|-y_1\cdot x_2)^{2\beta}}{|x_1||x_2|}\delta\binom{|x_1|+|x_2|-|y_1|-|y_2|}{x_1+x_2-y_1-y_2}\,\mathrm{d}x_1\mathrm{d}x_2,
\end{equation}
for $y=(y_1,y_2)\in\R^{2d}$; by changing variables $(x_2,y_2) \mapsto -(x_2,y_2)$, to complete the proof of the sharp inequality it suffices to prove the following key lemma.
\begin{lemma}\label{L1}
For $y\in\R^{2d}$ and $\beta>\frac{1-d}{4}$ we have that
\[
\mathrm{I}_\beta(y) = (2\pi)^{\frac{d-1}{2}} \frac{\Gamma(\frac{d-1}{2} + 2\beta)}{\Gamma(d-1+2\beta)}  (|y_1||y_2|-y_1\cdot y_2)^{\frac{d-3}{2}+2\beta}.
\]
\end{lemma}
\begin{remark} In the case $\beta=0$, Lemma \ref{L1} was proved for $d=3$ by Foschi in \cite{Foschi} and for general dimensions in \cite{BR}. In each case the proof proceeds by  applying an appropriate Lorentz transformation to reduce to the case $y_1=-y_2$; we shall see that in fact this line of argument allows us to obtain the result for more general values of $\beta$. We also remark that in the case $\beta=0$, Lemma \ref{L1} may be seen to hold via a direct calculation using the homogeneity of the delta measure, and is contained in Lemma 4.1 of the earlier paper \cite{FK}.
\end{remark}

\proof[Proof of Lemma \ref{L1}]
In order to shorten the formulas, in what follows we define $\tau=|y_1|+|y_2|$ and $\xi=y_1+y_2$. It is then easy to see that $\mathrm{I}_\beta(y)$ equals
\begin{align*}
& \int_{\R^{2d+2}}\frac{\delta(\sigma_1-|x_1|)}{|x_1|}\frac{\delta(\sigma_2-|x_2|)}{|x_2|}\delta\binom{\sigma_1+\sigma_2-\tau}{x_1+x_2-\xi}(|y_1||x_2|-y_1\cdot x_2)^{2\beta} \,\mathrm{d}\sigma_1\mathrm{d}\sigma_2\mathrm{d}x_1\mathrm{d}x_2 \\
& =\int_{\R^{d+1}}\int_{\R^{d+1}}\frac{\delta(\sigma_1-|x_1|)}{|x_1|}\frac{\delta(\sigma_2-|x_2|)}{|x_2|}\delta\left(\begin{pmatrix}\sigma_1 \\ x_1\end{pmatrix}+\binom{\sigma_2}{x_2}-\binom{\tau}{\xi}\right) \\
& \qquad\qquad\qquad\qquad\times(|y_1|\sigma_2-y_1\cdot x_2)^{2\beta} \,\mathrm{d}\binom{\sigma_1}{x_1}\mathrm{d}\binom{\sigma_2}{x_2}.
\end{align*}
We now introduce the Lorentz transformation $L$, given by
\begin{equation}\label{L-def}
L\binom{t}{x}=\binom{\gamma(t-v\cdot x)}{x+\big(\frac{\gamma-1}{|v|^2}v\cdot x-\gamma t\big)v}, \qquad (t,x) \in \R \times \R^d
\end{equation}
with $v:=-\frac{\xi}{\tau}$ and 
\begin{equation}\label{gamma}
\gamma:=\frac{1}{(1-|v|^2)^{\frac{1}{2}}}=\frac{\tau}{(\tau^2-|\xi|^2)^{\frac{1}{2}}}.
\end{equation}
It is not hard to check that
\[L\binom{(\tau^2-|\xi|^2)^{\frac{1}{2}}}{0}=\binom{\tau}{\xi},\]
that $|\op{det}L|=1$ and that the measure $|x|^{-1}\delta(t-|x|)$ is invariant under the transformation $L$. Applying the change of variables $\binom{\widetilde{\sigma_j}}{\widetilde{x_j}}=L^{-1}\binom{\sigma_j}{x_j}$ for $j=1,2$ it follows that
\begin{align*}
\mathrm{I}_\beta(y) & =\int_{\R^{2(d+1)}}\frac{\delta(\widetilde{\sigma_1}-|\widetilde{x_1}|)}{|\widetilde{x_1}|}\frac{\delta(\widetilde{\sigma_2}-|\widetilde{x_2}|)}{|\widetilde{x_2}|}\delta\left(L\begin{pmatrix}\widetilde{\sigma_1} \\ \widetilde{x_1}\end{pmatrix}+L\binom{\widetilde{\sigma_2}}{\widetilde{x_2}}-L\binom{(\tau^2-|\xi|^2)^{\frac{1}{2}}}{0}\right)\\
& \qquad\qquad\qquad\times\left(\binom{|y_1|}{-y_1}\cdot L\binom{\widetilde{\sigma_2}}{\widetilde{x_2}}\right)^{2\beta} \,\mathrm{d}\binom{\widetilde{\sigma_1}}{\widetilde{x_1}}\mathrm{d}\binom{\widetilde{\sigma_2}}{\widetilde{x_2}} \\
& = \int_{\R^{2d+2}}\frac{\delta(\widetilde{\sigma_1}-|\widetilde{x_1}|)}{|\widetilde{x_1}|}\frac{\delta(\widetilde{\sigma_2}-|\widetilde{x_2}|)}{|\widetilde{x_2}|}\delta\binom{\widetilde{\sigma_1}+\widetilde{\sigma_2}-(\tau^2-|\xi|^2)^{\frac{1}{2}}}{\widetilde{x_1}+\widetilde{x_2}} \\
& \qquad\qquad\qquad\times \left(\binom{|y_1|}{-y_1}\cdot L\binom{\widetilde{\sigma_2}}{\widetilde{x_2}}\right)^{2\beta}\,\mathrm{d}\widetilde{\sigma_1}\mathrm{d}\widetilde{\sigma_2}\mathrm{d}\widetilde{x_1}\mathrm{d}\widetilde{x_2} \\
& =\int_{\R^d}\frac{1}{|x|^2}\delta\left(2|x|-(\tau^2-|\xi|^2)^{\frac{1}{2}}\right)\left(\binom{|y_1|}{-y_1}\cdot L\binom{|x|}{x}\right)^{2\beta}\,\mathrm{d}x,
\end{align*}
where the last line follows by evaluating the integrals in $\widetilde{\sigma_1},\widetilde{\sigma_2}$ and $\widetilde{x_1}$ and then relabeling $\widetilde{x_2}=x$. We are now required to compute the quantity 
\[\binom{|y_1|}{-y_1}\cdot L\binom{|x|}{x}=|y_1|\gamma(|x|-v\cdot x)-y_1\cdot\left(x+\left(\frac{\gamma-1}{|v|^2}v\cdot x-\gamma |x|\right)v\right)\]
for $v=-\frac{\xi}{\tau}$, $\gamma$ given by \eqref{gamma} and $x$ on the support of the remaining delta measure; this is contained in the following.

\begin{lemma}\label{L2}
For each $y_1,y_2\in\R^d$ there exists $\omega_\ast\in\mathbb{S}^{d-1}$ such that
\[\binom{|y_1|}{-y_1}\cdot L\binom{|x|}{x}=\frac{|y_1||y_2|-y_1\cdot y_2}{2}\left(1+\frac{x}{|x|}\cdot \omega_\ast\right)\]
for any $x\in\R^d$ with $2|x|=(\tau^2-|\xi|^2)^{\frac{1}{2}}$.
\end{lemma}
Assuming Lemma \ref{L2} to be true for the moment, using polar co-ordinates, the definitions of $\tau$ and $\xi$, and then rotation invariance to take $\omega_\ast$ to $\mathbf{e}_1$,
\begin{align*}
\mathrm{I}_\beta(y) 
 = \frac{\left(|y_1||y_2|-y_1\cdot y_2\right)^{\frac{d-3}{2} + 2\beta}}{2^{2\beta+\frac{d-1}{2}}}\int_{\mathbb{S}^{d-1}}(1+\omega \cdot \mathbf{e}_1)^{2\beta}\,\mathrm{d}\omega.
\end{align*}
By changing variables, it follows that if $\beta>\frac{1-d}{4}$ then the integral term is
\begin{align*}
 |\mathbb{S}^{d-2}| \int_{-1}^1 (1+t)^{\frac{d-3}{2} + 2\beta}(1-t)^{\frac{d-3}{2}}\,\mathrm{d}t & = |\mathbb{S}^{d-2}|2^{d-2 + 2\beta} \mathrm{B}\left(\frac{d-1}{2} + 2\beta,\frac{d-1}{2}\right) \\
 & = 2^{d-1 + 2\beta} \pi^{\frac{d-1}{2}} \frac{\Gamma(\frac{d-1}{2} + 2\beta)}{\Gamma(d-1+2\beta)},
\end{align*}
where $\mathrm{B}$ is the classical beta function, and we have used the identity $\mathrm{B}(x,y) = \frac{\Gamma(x)\Gamma(y)}{\Gamma(x+y)}$ for appropriate $x$ and $y$. The claimed expression for $\mathrm{I}_\beta$ in the statement of Lemma \ref{L1} follows.  \endproof

It only remains to prove Lemma \ref{L2}.
\proof[Proof of Lemma \ref{L2}] After some straightforward calculations and simplifications we deduce that
\begin{align*}
L\binom{|x|}{x} & =\frac{1}{(\tau^2-|\xi|^2)^{\frac{1}{2}}}\binom{\tau|x|+\xi\cdot x}{x(\tau^2-|\xi|^2)^{\frac{1}{2}}+\xi\big(|x|+\xi\cdot x(\tau+(\tau^2-|\xi|^2)^{\frac{1}{2}})^{-1}\big)} \\
& = \frac{|x|}{(\tau^2-|\xi|^2)^{\frac{1}{2}}}\binom{\tau+\xi\cdot x^\prime}{x^\prime(\tau^2-|\xi|^2)^{\frac{1}{2}}+\xi\big(1+\xi\cdot x^\prime(\tau+(\tau^2-|\xi|^2)^{\frac{1}{2}})^{-1}\big)},
\end{align*}
where we recall $x^\prime:=\frac{x}{|x|}$. By our assumption $|x|=\frac{1}{2}(\tau^2-|\xi|^2)^{\frac{1}{2}}$,
\[
\binom{|y_1|}{-y_1}\cdot L\binom{|x|}{x}=\frac{1}{2}\left(|y_1|(\tau+\xi\cdot x^\prime) -2y_1\cdot x-y_1\cdot\xi\bigg(1+\frac{\xi\cdot x^\prime}{\tau+2|x|}\bigg) \right).
\]
Moreover, using the definitions $\tau=|y_1|+|y_2|$ and $\xi=y_1+y_2$,
\[|y_1|\tau-y_1\cdot\xi=|y_1||y_2|-y_1\cdot y_2,\]
and
\begin{align*}
& |y_1|\xi\cdot x^\prime-  2y_1\cdot x-\xi\cdot x^\prime\frac{\xi\cdot y_1}{\tau+2|x|} \\
& =x^\prime\cdot\left(|y_1|\xi-2|x|y_1-\xi\frac{\xi\cdot y_1}{\tau+2|x|}\right)  \\
& =\frac{1}{\tau+2|x|} x^\prime\cdot\bigg(\xi\big(|y_1|(\tau+2|x|)-\xi\cdot y_1\big)-2|x|(\tau+2|x|)y_1\bigg) \\
& =x^\prime\cdot\left(y_2\left(\frac{|y_1||y_2|-y_1\cdot y_2+2|x||y_1|}{|y_1|+|y_2|+2|x|}\right)-y_1\left(\frac{|y_1||y_2|-y_1\cdot y_2+2|x||y_2|}{|y_1|+|y_2|+2|x|}\right)\right)\\
& = x^\prime \cdot \frac{2|x|\big(y_2(|x|+|y_1|)-y_1(|x|+|y_2|)\big)}{|y_1|+|y_2|+2|x|}, 
\end{align*}
we get
\[
\binom{|y_1|}{-y_1}\cdot L\binom{|x|}{x}=\frac{1}{2}\left(|y_1||y_2|-y_1\cdot y_2+x^\prime\cdot z\right),
\]
where
\[
z=z(y_1,y_2):=\frac{2|x|\big(y_2(|x|+|y_1|)-y_1(|x|+|y_2|)\big)}{|y_1|+|y_2|+2|x|}.
\]
For $x$ such that  
\begin{equation} \label{e:xcondition}
2|x|^2 = \frac{1}{2}(\tau^2-|\xi|^2) =  |y_1||y_2|-y_1\cdot y_2,
\end{equation} 
we claim that 
\begin{equation} \label{e:zlength}
|z| = |y_1||y_2|-y_1\cdot y_2
\end{equation}
which would complete the proof of Lemma \ref{L2} by taking $\omega_\ast = z'$.

To see \eqref{e:zlength}, first note that by expanding the square,
\begin{align*}
& \left|2|x|\big(y_2(|x|+|y_1|)-y_1(|x|+|y_2|)\big)\right|^2 \\
& = 4|x|^2\left( |x|^2(|y_1|^2 + |y_2|^2 - 2y_1 \cdot y_2) + 2(|y_1||y_2| - y_1 \cdot y_2)(|y_1||y_2| + |x||y_1| + |x||y_2|) \right) 
\end{align*}
which may be rewritten as
\begin{align*}
& 2|x|^2\big(2|x|^2(|y_1|^2 + |y_2|^2 - 2y_1 \cdot y_2) \\
& \qquad + (|y_1||y_2| - y_1 \cdot y_2)(2|y_1||y_2| + 4|x||y_1| + 4|x||y_2|) + 4|x|^2|y_1||y_2|  \big). 
\end{align*}
Using \eqref{e:xcondition} repeatedly, this can be expressed as
\begin{align*}
& 2|x|^2\big(2|x|^2(|y_1|^2 + |y_2|^2) + 4|x|^2( |y_1||y_2| - y_1 \cdot y_2) \\
& \qquad + (|y_1||y_2| - y_1 \cdot y_2)(2|y_1||y_2| + 4|x||y_1| + 4|x||y_2|)  \big) 
\end{align*}
which is equal to
\begin{align*}
2|x|^2(|y_1||y_2| - y_1 \cdot y_2)(|y_1| + |y_2| + 2|x|)^2 = (|y_1||y_2| - y_1 \cdot y_2)^2(|y_1| + |y_2| + 2|x|)^2,
\end{align*}
and from this we obtain \eqref{e:zlength}.
\endproof

\subsection{Cases of equality in Theorem \ref{thm1}} 
The remaining task in the proof of Theorem \ref{thm1} is to establish the claim regarding the characterisation of the extremisers. 
As we have already observed, since the only place an inequality was used in the proof of \eqref{b-eq} was in the application of the arithmetic-geometric mean inequality, it follows that equality holds in \eqref{b-eq} if and only if (replacing $x_2$ and $y_2$ with $-x_2$ and $-y_2$ in \eqref{f-eq}) 
\[|x_1||x_2|\widehat{f}(x_1)\widehat{g}(x_2)=|y_1||y_2|\widehat{f}(y_1)\widehat{g}(y_2)\]
for almost every $(x_1,x_2,y_1,y_2)\in\R^{4d}$ satisfying 
$$
x_1+x_2 = y_1+y_2 \qquad \mbox{and} \qquad |x_1|+|x_2|=|y_1|+|y_2|. 
$$
Equivalently, equality holds if and only if
\begin{equation}\label{cone-eq}
h_1(x_1)h_2(x_2)=\Lambda(|x_1|+|x_2|,x_1+x_2)
\end{equation}
for almost every $x_1,x_2\in\R^d$ and some scalar function $\Lambda$, where $h_1:=|\cdot|\widehat{f}$ and $h_2:=|\cdot|\widehat{g}$, and, by symmetry and re-scaling, we can assume that $h_1=h_2=h$. 

The functional equation \eqref{cone-eq} was solved by Foschi in \cite{Foschi} in the case $d=3$ under the assumption that $h$ is locally integrable, and this argument was generalised to obtain the solution for all $d \geq 2$ (under the same local integrability assumption) in \cite{BR}. In particular, the set of locally integrable functions $h$ satisfying \eqref{cone-eq} have the form
\[
h(x) = e^{a|x| + b \cdot x + c}
\]
for some $a, c \in \mathbb{C}$ and $b \in \mathbb{C}^d$. 

Our goal will be to show that $|h|^\gamma$ is locally integrable for some $\gamma > 0$. Once established, it follows from \eqref{cone-eq} that
\[
|h(x)|^\gamma =  e^{a|x| + b \cdot x + c}
\]
for some $a, c \in \mathbb{R}$ and $b \in \mathbb{R}^d$, and taking $\gamma^{-1}$ powers, we obtain the desired form for $|h(x)|$. It then follows that $\frac{h}{|h|}$ is a bounded (hence locally integrable) solution of \eqref{cone-eq} and therefore
\[
\frac{h(x)}{|h(x)|} =  e^{ia'|x| + ib' \cdot x + ic'}
\]
for some $a', c' \in \mathbb{R}$ and $b' \in \mathbb{R}^d$. This gives the claimed characterisation of extremisers in Theorem \ref{thm1}.

To justify the remaining claim that $|h|^\gamma$ is locally integrable for some $\gamma > 0$, proceeding initially as in \cite{BR}, if $B$ is the euclidean ball centered at the origin of radius $N$, then by the Cauchy--Schwarz inequality, 
\begin{align*}
\left(\int_B |h| \right)^4 & \leq \mathcal{I}_\beta(f,f) \int_B\int_B(|x_1||x_2|)^{\frac{5-d}{2}-2\beta}(1-x_1^\prime\cdot x_2^\prime)^{\frac{3-d}{2}-2\beta}\,\mathrm{d}x_1\mathrm{d}x_2 \\
& =  \mathcal{I}_\beta(f,f) \left(\int_{\mathbb{S}^{d-1}}\int_{\mathbb{S}^{d-1}}(1-\omega_1\cdot\omega_2)^{\frac{3-d}{2}-2\beta}\,\mathrm{d}\omega_1\mathrm{d}\omega_2\right)\left(\int_0^N r^{\frac{d+3}{2}-2\beta}\,\mathrm{d}r\right)^2.
\end{align*}
The integral in $r$ is finite when $\beta<\frac{d+5}{4}$, but the integral over the sphere is equal 
\[
C_d\int_{-1}^1 (1-t)^{\frac{d-3}{2}}(1+t)^{-2\beta}\,\mathrm{d}t,
\]
which is finite for $\beta<\frac{1}{2}$. 

Now suppose $\beta \geq \frac{1}{2}$, for which we make use of the reverse Hardy--Littlewood--Sobolev inequality on the sphere (for a proof see \cite{Beckner}), which states that if $g\geq 0$ then
\begin{equation*}
H_\lambda(g) \gtrsim \|g\|_{L^p(\mathbb{S}^{d-1})}^2,
\end{equation*}
where $\lambda<0$ and $p:=\frac{2(d-1)}{2(d-1)-\lambda}\in\left(0,1\right)$. Using \eqref{HLStrick}, this means
\[
\mathcal{I}_\beta(f,f) \geq C_{\beta,d} \|T_\beta f\|_{L^p(\mathbb{S}^{d-1})}^2
\]
where $p = \frac{2(d-1)}{3d-5+4\beta} \in (0,\frac{2}{3}]$, and thus it suffices to show
\begin{equation} \label{e:lastone}
\int_B |h|^{2p} \leq C_{\beta,d,N} \|T_\beta f\|_{L^p(\mathbb{S}^{d-1})}^p.
\end{equation}
Applying H\"older's inequality and Minkowski's inequality for integrals, using that $\frac{1}{p} > 1$, we obtain
\begin{align*}
\int_B |h(x)|^{2p} \,\mathrm{d}x 
& = \int_0^N \left(\int_{\mathbb{S}^{d-1}}|\widehat{f}(r\omega)|^{2p} r^{2p+d-1} \,\mathrm{d}\omega\right)\,\mathrm{d}r \\
& \leq C_{\beta,d,N}\left(\int_0^N\left(\int_{\mathbb{S}^{d-1}}|\widehat{f}(r\omega)|^{2p}r^{2p+d-1}\,\mathrm{d}\omega\right)^{\frac{1}{p}} \,\mathrm{d}r\right)^{p}\\
& \leq C_{\beta,d,N} \int_{\mathbb{S}^{d-1}}\left(\int_0^N |\widehat{f}(r\omega)|^2 r^{2 + \frac{d-1}{p}}\,\mathrm{d}r\right)^{p}\mathrm{d}\omega.
\end{align*}
Since $2 + \frac{d-1}{p} > \frac{3}{2}(d-1) + 2\beta$, it follows that $r^{2 + \frac{d-1}{p}} \leq C_N r^{\frac{3}{2}(d-1) + 2\beta}$ for all $r \in [0,N]$, and \eqref{e:lastone} follows.

\endproof

\section{Sharp bilinear estimates for ($++$) waves}

In this last section, we present analogous results to Theorem \ref{thm1}, Corollary \ref{c:radial} and Corollary \ref{cor4} for $(++)$ waves instead of $(+-)$ waves.
\begin{theorem}\label{t:plusplus}
Let $d\geq 2$ and $\beta \in \mathbb{R}$. Then,
\begin{equation}\label{b-eq2}
\||\square|^\beta(e^{itD}fe^{itD}g)\|_{L^2(\R^{d+1})}^2\leq \mathbf{W}^\prime(\beta,d)\mathcal{I}_\beta(f,g)
\end{equation}
holds with constant
\[
\mathbf{W}^\prime(\beta,d) = \frac{2^{\frac{5-7d + 4\beta}{2}}\pi^{\frac{2-5d}{2}}}{\Gamma(\frac{d}{2})}.
\]
For $\beta >\frac{2-d}{2}$ the constant is sharp and equality holds if and only if
$f=\lambda g$ satisfy \eqref{maxdef} for some $\lambda, a, c\in\mathbb{C}$ with $\operatorname{Re}a<0$, and $b\in\mathbb{C}^d$ with $|\operatorname{Re}b|<-\operatorname{Re}a$.
\end{theorem}

Note that the range of $\beta$ for the sharp estimates in Theorem \ref{t:plusplus} is the same as in Theorem \ref{thm1} when $d=2$, but is strictly larger otherwise. This is due to the presence of additional symmetry in the left-hand side of \eqref{b-eq2} which does not appear to occur in \eqref{b-eq}, permitting a simpler argument for which the restriction in Theorem \ref{thm1} does not arise. This is reminiscent of recent work in \cite{BBJP} where bilinear estimates (with and without complex conjugate) are studied for the Schr\"odinger evolution operator $e^{it\Delta}$. In that case, however, the presence of a complex conjugate causes a change in the shape of the multiplier at the level of sharp estimates and moreover the class of extremisers for the estimates analogous to \eqref{b-eq} and \eqref{b-eq2} are shown to be different. 

Since the right-hand side of \eqref{b-eq2} is in terms of the same functional $\mathcal{I}_\beta(f,g)$, the arguments in Section \ref{s:cor} immediately give the following.
\begin{corollary} 
Let $d \geq 2$, $\beta > \frac{2-d}{2}$ and let
\[
\mathbf{C}'(\beta,d) = 2^{4\beta - 1} \pi^{\frac{1-d}{2}} \frac{\Gamma(d-2 + 2\beta)}{\Gamma(\frac{3d-5}{2} + 2\beta)}.
\]
\emph{(i).} If $f$ and $g$ are radially symmetric, then 
\begin{equation*}
\left\||\square|^\beta (e^{itD}f e^{itD}g)\right\|_{L^2(\R^{d+1})}^2\leq \mathbf{C}'(\beta,d)\|f\|_{\dot{H}^{\frac{d-1}{4} + \beta}(\R^d)}^2\|g\|_{\dot{H}^{\frac{d-1}{4} + \beta}(\R^d)}^2
\end{equation*}
where the constant is sharp, with equality if and only if $f=\lambda g$ satisfy \eqref{maxdef} with $\lambda, c\in\mathbb{C}$, and $a \in \mathbb{C}$ such that $\operatorname{Re}a < 0$, and $b=0$.

\emph{(ii).} If $\beta \in (\frac{3-d}{4},\frac{5-d}{4}]$, then 
\begin{equation*}
\left\||\square|^\beta(e^{itD}f)^2\right\|_{L^2(\R^{d+1})}^2\leq \mathbf{C}'(\beta,d)\|f\|_{\dot{H}^{\frac{d-1}{4} + \beta}(\R^d)}^4
\end{equation*}
and the constant is sharp, with equality if and only if $f$ satisfies \eqref{maxdef} with $c\in\mathbb{C}$, and $a \in \mathbb{C}, b \in \mathbb{C}^d$ are such that $\operatorname{Re}a < 0$ and $\operatorname{Re}b=0$.

\emph{(iii).} If $\beta \in (\frac{2-d}{4},\frac{3-d}{4}]$ and $p=\frac{2(d-1)}{3d-5+4\beta}$, then 
\begin{equation*}
\left\||\square|^\beta (e^{itD}f e^{itD}g)\right\|_{L^2(\R^{d+1})}^2\leq \mathbf{C}'(\beta,d)|\mathbb{S}^{d-1}|^{\frac{3-d-4\beta}{d-1}}\|T_\beta f\|_{L^p}\|T_\beta g\|_{L^p}
\end{equation*}
and the constant is sharp, with equality if and only if $f=\lambda g$ satisfy \eqref{maxdef} with $\lambda, a, c\in\mathbb{C}$ and $b\in\mathbb{C}^d$ where $|\operatorname{Re}b|<-\operatorname{Re}a$. 
\end{corollary}

We end by briefly indicating how to modify the argument in Section \ref{s:proofs} to prove Theorem \ref{t:plusplus}. As noted earlier, it is also possible to follow the earlier arguments in \cite{Foschi} and \cite{BR}, but the advantage of our unified approach is that it highlights the additional difficulty in the case of $(+-)$ waves.

Firstly,
\begin{align*}
& 2^{-2\beta}(2\pi)^{3d-1}\||\square|^\beta (e^{itD}fe^{itD}g)\|_{L^2(\R^{d+1})}^2 \\
& = \int_{\R^{4d}}(|y_1||y_2|-y_1\cdot y_2)^{2\beta}\widehat{f}(x_1)\widehat{g}(x_2)\overline{\widehat{f}(y_1)\widehat{g}(y_2)}\delta\binom{|x_1|+|x_2|-|y_1|-|y_2|}{x_1+x_2-y_1-y_2}\,\mathrm{d}x\mathrm{d}y \\
& = \int_{\R^{4d}}(|y_1||y_2|-y_1\cdot y_2)^{2\beta} \frac{\widehat{F}(x)\overline{\widehat{F}(y)}}{(|x_1||x_2||y_1||y_2|)^{\frac{1}{2}}}\delta\binom{|x_1|+|x_2|-|y_1|-|y_2|}{x_1+x_2-y_1-y_2}\,\mathrm{d}x\mathrm{d}y,
\end{align*}
but in this case we do not need to relabel any indices and we take $\widehat{F}:=(|\cdot|^{\frac{1}{2}}\widehat{f})\otimes(|\cdot|^{\frac{1}{2}}\widehat{g})$. Proceeding exactly as in the proof of \eqref{b-eq} we then obtain
\[
\frac{\||\square|^\beta (e^{itD}fe^{itD}g)\|_{L^2}^2}{2^{2\beta}(2\pi)^{1-3d}}\leq\int_{\R^{2d}}(|y_1||y_2|-y_1\cdot y_2)^{2\beta}|\widehat{F}(y)|^2 \mathrm{I}_0(y) \,\mathrm{d}y, 
\] 
where we have used the fact that $|x_1||x_2|-x_1\cdot x_2=|y_1||y_2|-y_1\cdot y_2$ on the support of the delta measures in this case, and $\mathrm{I}_0$ is defined in \eqref{e:Idefn}. Hence (for any admissible $\beta$) we need use only the $\beta=0$ case of Lemma \ref{L1} to obtain \eqref{b-eq2}. From a technical point of view, Lemma \ref{L1} is substantially easier in the case $\beta = 0$ and this explains why the $(++)$ case is easier.

Finally, the characterisation of extremisers follows in precisely the same manner, since the above argument shows that equality holds in $\eqref{b-eq2}$ if and only if
\[
|x_1||x_2|\widehat{f}(x_1)\widehat{g}(x_2)=|y_1||y_2|\widehat{f}(y_1)\widehat{g}(y_2)
\]
almost everywhere on the support of the delta measures, and this is the same functional equation which arose in the $(+-)$ case.

\begin{acknowledgement}
The first author was supported by JSPS Kakenhi grant number 26887008. The second author acknowledges support from the Universitas 21 PhD Scholarship 2014 (University of Birmingham). He also wishes to thank Waseda University, where part of this research was conducted, for their kind hospitality.
\end{acknowledgement}

\end{document}